\documentclass[reqno,10pt]{amsart}
\usepackage{amscd,amssymb,amsmath,amscd,misha,stmaryrd}
\usepackage{latexsym}
\usepackage[all]{xy}
\usepackage[colorlinks=true]{hyperref}

\newcommand{\Log}{{\mathbf{Log}}}
\def\Mon{{\rm Mon}}
\def\int{{\rm int}}
\def\sat{{\rm sat}}
\def\ket{{\rm ket}}
\def\rmlog{{\rm log}}
\def\..{,\dots,}

\begin{document}

\author{Michael Temkin}
\title{Introduction to logarithmic geometry}
\address{\tiny{Einstein Institute of Mathematics\\
               The Hebrew University of Jerusalem\\
                Edmond J. Safra Campus, Giv'at Ram, Jerusalem, 91904, Israel}}
\email{\scriptsize{michael.temkin@mail.huji.ac.il}}
\thanks{This research is supported by  BSF grants 2014365 and 2018193, ERC Consolidator Grant 770922 - BirNonArchGeom.}

\setcounter{tocdepth}{1}

\maketitle

\section{Introduction}
These notes will substitute a chapter in a book on recent advances in resolution of singularities based on a series of minicourses given at an Oberwolfach seminar. It is based on a minicourse on logarithmic resolution of singularities given by the author, and it provides an extended version of its first part devoted to introduction to logarithmic geometry with a view towards applications to resolution. I do not aim to build a theory with proofs (and this is impossible in a 3-4 lecture long course). The goal is to make the reader familiar with basic definitions, constructions, techniques and results of logarithmic geometry. I formulate most of the results as ``Exercises'' and try to keep them at a reasonable level of difficulty. References to the literature are also provided. At the first reading of the material it may be worth just to read the formulations and hints or comments about main ideas of the arguments, without trying to solve them or read proofs in the cited papers.

\subsection{History and motivation}

\subsubsection{The discovery}
Logarithmic structures and schemes were discovered by J.-M. Fontaine and L. Illusie on Sunday, July 17, 1988 during a discussion in a train on their travel to Oberwolfach workshop "Aritmetische Algebraische Geometrie". In fact, the discussion was in the continuation of an IHES seminar that had taken place in the spring, and the construction was motivated by the necessity of finding a suitable framework in which an analogue of Steenbrink's limiting Hodge structure for a semistable reduction over a complex disc could be defined in mixed characteristic in order to make sense of the $C_{\rm st}$-conjecture of Fontaine and Jannsen.

During the workshop Illusie prepared a short summary of the discussion and showed it to K. Kato, who was very enthusiastic about the new notion and very quickly wrote the first paper, where these notions were introduced: ``Logarithmic structures of Fontaine-Illusie''. The new theory turned out to be extremely useful because of the following features:

\begin{itemize}
\item[(1)] It provides a more general notion of smoothness, which allows to work with many classical non-smooth objects similarly to the smooth ones. In particular, it conceptually adjusts various cohomology theories to this generalized context.
\item[(2)] It conceptually treats various notions of boundaries, such as normal crossings divisors, and it often provides a functorial way to compactify various moduli spaces -- smooth objects often degenerate to logarithmically smooth (but non-smooth) objects over the boundary.
\item[(3)] It provides a conceptual way to bookkeep information on closed subspaces and fibers, in particular, leading to better solutions of deformation problems.
\end{itemize}
Off course these three large classes of properties are tightly connected and often show up altogether.

\subsubsection{Precursors}
In fact, log geometry had numerous precursors, which it absorbed and generalized. Without pretending to provide a full list, here are a few most important ones, which will be discussed in \S\ref{precursors} in more detail:
\begin{itemize}
\item[(1)] Normal crossings divisors, especially, when viewed as a boundary used to compactify a smooth variety, correspond to log structures. In fact, a smooth variety with a normal crossing divisor is nothing else but a log smooth variety, which is smooth.
\item[(2)] Toroidal geometry, which was introduced in \cite{KKMS} to prove semistable reduction theorem, is, in fact, the theory of log smooth log varieties. Toroidal morphisms between toroidal varieties are nothing else but log smooth morphisms.
\item[(3)] Deligne's generalized divisors and logarithmic structures of Deligne-Faltings, \cite[Complement 1]{Kato-log}. See also \cite[Section 2]{Deligne_garden}, where Deligne's letter and its influence on the discovery of log schemes is described.
\item[(4)] Logarithmic differentials, logarithmic versions of various complexes, etc., which were defined ad hoc, obtain a conceptual interpretation in log geometry.
\item[(5)] Semistable morphisms are log smooth, so semistable reduction theorem literally becomes a desingularization theorem in log geometry. Furthermore, the snc divisor sitting in the closed fiber of a semistable family is log smooth over the log point -- a more exotic object, which bookkeeps the log structure purely algebraically. In a sense, a log point is a logarithmic analogue of usual non-reduced (or fat) points in the theory of schemes.
\end{itemize}

In fact, the role of log geometry in some other classical problems is being gradually clarified even nowadays. For example, it was clear that it is involved in resolution of singularities, at least through the exceptional divisor, but a careful study of this question not only shed a new light on known methods, but also led to discovery of a new generation of methods, which will be discussed in another chapter. Also, I think that the role of log geometry in compactifying various moduli spaces is not fully exploited yet, and it will increase in future research. This is tangential to the material we want to cover, so we will only discuss a couple of such examples in the sequel and refer the reader to \cite{logmoduli}, where the theory of log schemes is described from the point of view of applications to the theory of moduli spaces. Our main motivation is application to resolution of singularities in the chapter on logarithmic and relative resolution of singularities.

\subsection{Structure of the chapter}

\subsubsection{Overview}
The chapter starts with Section \ref{precursors} on precursors of logarithmic geometry: we discuss the log structure encoded by snc divisors, which is probably the first time log structures implicitly showed up in mathematics, and then recall the more general theory of toroidal varieties. Log schemes are introduced in \S\ref{logstrsec}. We first study necessary properties of monoids and introduce log structures, and in the end of the section we discuss log regular log schemes. Various properties of morphisms: charts, log smoothness and log \'etaleness, log differentials and log blowings up are reviewed in \S\ref{morsec}. Finally, in Section~\ref{Logsec} we discuss Olsson's stacks $\Log_X$ and the technique of reducing log geometry to geometry of stacks. This is the most technically demanding section, and the only one in which stacks are used. It will be used in the construction of the relative desingularization functor, but not in the more basic case of the absolute logarithmic desingularization.

\subsubsection{References and sources}
The single reference with all foundations worked out in detail is the recent book of Arthur Ogus \cite{logbook}. Originally, logarithmic geometry was established by Kazuya Kato in \cite{Kato-log}, and the theory of log regular log schemes, their desingularization and log blowings up was developed in \cite{Kato-toric} and \cite{Niziol}. Stacks $\Log_X$ were introduced and their relation to logarithmic properties was studied by Martin Olsson in \cite{Olsson-logarithmic}, and some further results were obtained in \cite{Molcho-Temkin}.

\subsubsection{Conventions}
Often we write ``log'' instead of ``logarithmic''.

\subsubsection{Acknowledgments}
I am very grateful to Luc Illusie for telling the story of discovery of log schemes and for reading the notes and making many helpful comments.

\tableofcontents

\section{Precursors}\label{precursors}
In this section we discuss some situations, where log structures are implicit actors. Later on they will serve as a source of examples and illustrations.

\subsection{Normal crossings divisors}

\begin{defin}
Let $X$ be a regular scheme and $D\into X$ a divisor.

(i) One says that $D$ is {\em strictly (or simple) normal crossings} or just {\em snc} at $x\in X$ if locally one has that $D=V(t_1\ldots t_s)$, where $t_1\..t_n$ is a regular family of parameters at $x$ and we use notation $V(f_1\..f_n)=\Spec_X(\calO_X/(f_1\..f_n))$.

(ii) One says that $D$ is {\em normal crossings} at $x$ if it is snc \'etale-locally at $x$.

A divisor is called {\em normal crossings} or {\em snc} if this is so everywhere on $X$.
\end{defin}

\begin{exer}\label{sncexer}
(i) Show that the number $s=s(x)$ in the definition is the number of branches of $D$ at $s$ (e.g. the number of irreducible formal components) and hence is an invariant of $D$ at $x$, also called the {\em multiplicity} of $D$ at $x$.

(ii) Define a stratification of $D$ by the multiplicity and show that each stratum $D(s)$ is a regular locally closed subscheme.

(iii) Assume that $X$ is regular and $D=\cup_{i\in I} D_i$ is a reduced divisor with irreducible components $D_i$. Show that $D$ is snc if and only if for each $J\subseteq I$ the scheme theoretic intersection $D_J=\cap_{j\in J}D_j$ is regular.
\end{exer}

\begin{rem}
Local computations with nc (resp. snc) divisors are done using \'etale (resp. Zariski) neighborhood, where it is given by the vanishing of the product of a subset of a family of regular parameters.
\end{rem}

\subsection{Toroidal schemes}

\subsubsection{Toric schemes}
Let $M$ be a lattice and $N=\Hom(M,\bbZ)$ the dual lattice. This can be encoded by the non-degenerate pairing $M\times N\to\bbZ$, but despite the symmetry, we will always view $N$ and $N_\bbR=N\otimes\bbR$ as geometric spaces, while the elements of $M$ will be viewed as functions $N\to\bbZ$ or $N_\bbR\to\bbR$.

Let $\sigma\in N_\bbR$ be an $M$-rational polyhedral cone, i.e. a cone given by finitely many conditions $m_i(x)\ge 0$ with $m_i\in M$. One associates to $\sigma$ an affine toric $k$-variety $\bfT_\sigma$ as follows. The dual cone $\sigma^\vee=\{z\in M_\bbR|\ z(\sigma)\ge 0\}$ is also rational and it follows easily that the monoid $M_\sigma:=\sigma^\vee\cap M$ is finitely generated. Therefore, $A_\sigma:=k[M_\sigma]$ is an affine $k$-algebra and $\bfT_\sigma:=\Spec(A_\sigma)$ is an affine variety. Note that $\bfT_\sigma$ is provided with the natural action of the torus $\bfT_0=\bfD_M=\Spec(k[M])$ whose lattice of characters is $M$. Furthermore, $\bfT_\sigma$ contains an open orbit isomorphic to $\bfT_0$, and this is the source of the terminology. The information encoded in the pair $(M,\sigma)$, often called the combinatorial information, is equivalent to the information encoded in $\bfT_\sigma$ with the torus action, and there is a very tight and natural relation between the combinatorial and geometric pictures:

\begin{exer}
(i) Show that $(M,\sigma)$ can be reconstructed from the affine toric variety $\bfT_\sigma$ as follows: $M$ is the lattice of characters of $\bfT_0$ and giving the $\bfT_0$-action on $\Spec(A)$ is equivalent to providing an $M$-grading of $A$, where $m$-homogeneous elements are the equivariant ones with the character $m$: they are acted on via the rule $t(f)=m(t)f$; the cone $\sigma$ is determined by the monoid $M_\sigma$, which is precisely the set of characters $m\in M$ with a non-zero homogeneous component $(A_\sigma)_m$.

(ii) Show that the set of orbits of $\bfT_\sigma$ is in a natural one-to-one correspondence with the faces of $\sigma$. In particular, the open orbit $\bfT_0$ corresponds to the vertex $0$ of $\sigma$.
\end{exer}

\begin{rem}
The affine theory can be globalized as follows. On the schematic part of the picture one glues affine toric varieties $T_i$ with the same torus $\bfD_M$ along isomorphic open affine toric subvarieties. A resulting object is called a toric scheme with respect to the torus $D_M$. On the combinatorial side one glues polyhedral cones $\sigma_i\subset N_\bbR$ along faces. The resulting objects are called {\em cone complexes}. Sometimes one only considers complexes embedded in $N_\bbR$ and calls them {\em polyhedral fans}. If the union of all faces of a fan is the whole $N_\bbR$, the fan is called a subdivision. In fact, any separated toric variety corresponds to a fan and a proper toric variety corresponds to a subdivision. We do not go into details and recommend a standard literature, e.g. \cite{Fulton}.
\end{rem}

\subsubsection{Monoidal resolution of singularities}
In the geometry of polyhedral complexes regular simplicial cones play the role of non-singular points.

\begin{defex}
(i) A polyhedral cone $\sigma$ is called {\em regular} if its sharpening $\oM_\sigma=M_\sigma/(M_\sigma)^\times$ is a free monoid $\bbN^r$ and hence $M_\sigma\toisom\bbN^r\oplus\bbZ^s$. Show that this happens if and only if there exists a basis $e_1\..e_n$ of $N$ such that each edge of $\sigma$ contains some $e_i$. Thus, $\sigma$ is a simplicial cone of a very special form.

(ii) A polyhedral complex is called {\em regular} if all its cones are.
\end{defex}

A combinatorial or monoidal resolution of singularities is the following result:

\begin{theor}
Any polyhedral complex possesses a regular subdivision.
\end{theor}

\begin{rem}
We will not prove this theorem, but only make a couple of remarks. One can construct such subdivision by applying functorial (hence equivariant) resolution of singularities to toric varieties and translating it to the combinatorial language. Clearly, this is a too complicated solution that does not admit a simple combinatorial interpretation. A construction of a simple canonical solution of this problem was missing in the literature until very recently, see \cite[Theorem~4.6.1]{Jarek-torres}, but various non-canonical solutions were well known. First, the claim easily reduces to the case of a single cone. Second, using the barycentric subdivision one reduces to the case of a simplicial cone. Then one defines some invariants of the singularity (essentially, they measure the discrepancy between $\sigma\cap N$ and the monoid generated by the elements of $N$ lying on the edges) and finds simplicial subdivisions that decrease it.
\end{rem}

\subsubsection{Toroidal embeddings}
Toroidal varieties (or schemes) \'etale-locally (or formally locally) are modelled on toric schemes. This turns out to be sufficient to extend various toric constructions, such as toric blowings up and toric resolution of singularities, to much wider context. The miracle enabling this is that these constructions pull back to the same operation on the toroidal scheme independently of the choice of a toric chart.

For simplicity we will only consider the case of varieties. There is no torus action anymore, but it turns out that a large portion of the structure can be encoded just in the open orbit:

\begin{defin}
(i) A {\em toroidal variety} (or a toroidal embedding) over a field $k$ is a pair $(X,U)$ with $X$ a $k$-variety and $U\into X$ an open subscheme such that \'etale locally $X$ possesses an \'etale morphism to a toric scheme such that $U$ is the preimage of the torus. Namely, there exists an \'etale covering $\coprod_i X_i\to X$ and \'etale morphisms $h_i\:X_i\to\bfT_{\sigma_i}$, called {\em toroidal chart}, such that $U\times_XX_i$ is the preimage of the torus of $\bfT_{\sigma_i}$. If, moreover, the covering can be chosen to be Zariski, say $X=\cup_iX_i$, then the toroidal variety is called {\em simple} (or without self-intersections).

(ii) A morphism of toroidal varieties $(Y,V)\to(X,U)$ is any morphism $Y\to X$ taking $V$ to $U$.
\end{defin}

\begin{exam}
If $X$ is regular, $D$ is a normal crossings divisor and $U=X\setminus D$, then $(X,U)$ is a toroidal variety modelled on regular simplices $\sigma$, and the free monoid $M_\sigma$ is generated by the regular parameters which define the branches of $D$ (on an appropriate \'etale neighborhood).
\end{exam}

\subsubsection{(Non-)uniqueness of charts}
A very natural question is to what extent the charts are unique. We will show below that the toric monoid is essentially unique, and the chart is unique up to units. For simplicity, assume that $i\:U\into X$ is a simple toroidal variety and set $\calM=i_*\calO_U^\times\cap\calO_X$ and $\ocalM=\calM/\calO_X^\times$ (in the general case one would have to use \'etale sheaves, as we will do in the section about log schemes). The latter sheaf is the sheaf of {\em toroidal Cartier divisors}, i.e. divisors supported on $D=X\setminus U$.

\begin{exer}
(i) Show that any local chart that maps $x$ to the closed stratum of some $\bfT_\sigma$ induces in isomorphism $M_\sigma=\ocalM_x$, in particular, $\bfT_\sigma$ depends only on $(X,U)$ and $x$.

(ii) Conversely, any section $s\:\ocalM_x\to\calM_x\subset\calO_{X,x}$ of $\calM_x\onto\ocalM_x$ can be extended to a local chart on a neighborhood of $x$. (Hint: in addition to $s$ one should choose a regular family of parameters on the stratum of $D=X\setminus U$ through $x$ and lift them to elements of $\calO_{X,x}$.)
\end{exer}


\section{Logarithmic structures and schemes}\label{logstrsec}
In this section we introduce the category of log schemes and study its basic properties.

\subsection{Monoids}
Unless said to the contrary, by a monoid we always mean a commutative additively written monoid $M=(M,+,0)$. By $M^\times$ we denote the subgroup of invertible elements and the {\em sharpening} of $M$ is $\oM=M/M^\times$. One says that $M$ is {\em sharp} if $M=\oM$.

\subsubsection{Basic constructions}
All categories of algebraic objects, such as groups, rings, commutative rings, etc. are complete and cocomplete -- possess all small limits and colimits. Furthermore, limits are compatible with set-theoretic limits, and colimits are obtained using generators and relations. In particular, this is true for the category $\Mon$ of monoids. The main examples of limits and colimits we will use are as follows:

\begin{exer}\label{basicmonoids}
(i) $M\times N$ is just the usual product of sets with componentwise addition and it coincides with the coproduct, usually denoted $M\oplus N$.

(ii) For homomorphisms $f\:M\to L$, $g\:N\to L$ the fiber product $M\times_L N$ is the submonoid of $M\times N$ given by $f(m)=g(n)$.

(iii) Instead of kernels in the category of monoids one uses {\em congruence relations}, that is, equivalence relations $R\subseteq M\times M$ which are also submonoids: if $M\to N$ is a surjective homomorphism, then the induced equivalence relation $R\subseteq M\times M$ is a submonoid, and conversely any congruence relation appears in this way.

(iv) For homomorphisms $f\:L\to M$, $g\:L\to N$ the pushout $M\oplus_L N$ is the quotient of $M\oplus N$ by the minimal congruence relation such that $f(l)\sim g(l)$ for any $l\in L$.
\end{exer}

\subsubsection{Ideals}
An {\em ideal} $I\subseteq M$ is a subset such that $I+M=I$, where the convention is that $I=\emptyset$ is also an ideal (the analogue of the zero ideal of rings). An ideal $I$ is {\em prime} if $x+y\in I$ implies that $x\in I$ or $y\in I$. The set of all prime ideals is denoted $\Spec(M)$ and called the {\em fan} of $M$. An ideal of the form $(a)=a+M$ is called {\em principal}.

\begin{exer}
Show that taking the preimage establishes a bijection between the ideals of $M$ and $\oM$.
\end{exer}

\subsubsection{Fine monoids}
We will usually work with fine monoids:

\begin{defin}
(i) $M$ is finitely generated if there exists a surjective homomorphism $\bbN^l\to M$.

(ii) $M$ is {\em integral} (or cancellative) if $n+m=n+m'$ implies that $m=m'$.

(iii) $M$ is {\em fine} if it is integral and finitely generated.
\end{defin}

\subsubsection{The Grothendieck group}
Recall that there is a canonical way to turn monoid into a group:

\begin{defex}
(i) Show that there is a universal homomorphism $M\to M^\gp$ with $M^\gp$ a group, called the {\em Grothendieck group} of $M$. (Hint: for example, one can bound the cardinality of $M^\gp$ because it is generated by the image of $M$, and then general representability theorems do the job because the category of groups is complete.)

(ii) Construct $M^\gp$ explicitly as the quotient of $M^2$ by the following equivalence relation: $(m,n)\toisom(m',n')$ if there exists $l\in M$ such that $l+m+n'=l+m'+n$.

(iii) Show that $M$ is integral if and only if the homomorphism $M\to M^\gp$ is injective.

(iv) Define the {\em integralization} $M^\int$ to be the image of $M$ in $M^\gp$. Show that $M\to M^\int$ is the universal homomorphism from $M$ to an integral monoid and $M\mapsto M^\int$ is the left adjoint functor to the embedding of the category of integral monoids $\Mon^\int$ into $\Mon$.
\end{defex}

\subsubsection{Fs monoids}
An especially nice class of monoids is defined as follows:

\begin{defin}
(i) An integral monoid $M$ is {\em saturated} if for each $m\in M^\gp$ with $am\in M$ for $m\in\bbZ_{>0}$ one has that $m\in M$.

(ii) A fine saturated monoid is called {\em fs}.

(iii) A sharp fs monoid is called {\em toric}.
\end{defin}

\begin{exer}
(i) Show that to give a toric monoid $M$ is equivalent to give a lattice $M^\gp$ and a rational polyhedral cone $\sigma$ in $M^\gp\otimes\bbR$ such that $M=M^\gp\cap\sigma$.

(ii) Show that $\Spec(M)$ is naturally bijective to the set of faces of $\sigma$.
\end{exer}

\begin{exer}
Let $\Mon^\sat$ denote the category of saturated monoids. Show that the embedding $\Mon^\sat\into\Mon$ possesses a left adjoint functor, which is called the saturation functor and denoted $M\mapsto M^\sat$. Show that $M^\sat$ is just the saturation of $M^\int$ in $M^\gp$, that is the divisible hull of $M^\int$ in $M^\gp$.
\end{exer}

\subsection{Logarithmic structures}

\begin{defin}
Let $\tau$ be one of the following topologies -- Zariski, \'etale or flat.

(i) A $\tau$-prelogarithmic structure on a scheme $X$ is a sheaf of monoids $\calM$ on the site $X_\tau$ with a {\em structure homomorphism} of monoids $u:\calM\to (\calO_X,\cdot)$. A homomorphism of prelog structures $\calM\to \calM'$ is a homomorphism of sheaves of monoids compatible with the structure homomorphisms.

(ii) A $\tau$-logarithmic structure is a $\tau$-prelogarithmic structure which induces an isomorphism $u^{-1}(\calO_{X_\tau}^\times)\toisom\calO_{X_\tau}^\times$, and hence also $\calM^\times\toisom\calO_{X_\tau}^\times$. The sharpening $\ocalM=\calM/\calO_{X_\tau}^\times$ is called the {\em characteristic monoid} of $\calM$.

(iii) The default topology in this definition is the \'etale topology, so usually it will not be mentioned. A log structure $\calM$ induces a Zariski log structure $\calM_\Zar$ just by restricting. By a slight abuse of language we say that $\calM$ itself is {\em Zariski} if this restriction does not loose information, that is, $\calM=\veps^*(\calM_\Zar)$ for the morphism of sites $\veps\:X_\et\to X_\Zar$.
\end{defin}

\begin{rem}
(i) The homomorphism $u$ is an analog of exponentiation, and this is one of the reasons to use additive monoids as the source. Traditionally it is denoted $\alp$, say $\alp(m)=x$, but probably the exponential notation, such as $x=u^m$ instead of $u(m)$, is more suggestive. Informally, any such $m$ can be viewed as a branch of the logarithm of $x$, so the log structure can be viewed as fixing a monoid of branches of logarithms.

(ii) The \'etale topology is used instead of the Zariski topology first of all in order to adequately treat toroidal embeddings, which are not strict. As a rule, this might pose mild technical inconveniences, which can be bypassed. For example, Kato restricted the generality to Zariski log structures in \cite{Kato-toric}, but these results were generalized to the general case in \cite{Niziol}.

(iii) Fppf log structures are sometimes needed to bypass positive characteristic problems, usually by extracting appropriate $p$-th roots. They are rarely used and will not show up in these notes.
\end{rem}

To get an initial feeling let us consider some examples.

\begin{exam}\label{logstrexam}
(0) The minimal or {\em trivial log structure} is just $\calM=\calO_{X_\et}^\times$.

(1) The largest log structure with an injective $u$ is $\calM=\calO_{X_\et}$, but it is not too useful. The main exception is when $X$ is a ``small'' scheme -- a semi-local curve or the spectrum of a valuation ring, with the most useful case being when $X$ is a trait.

(2) A very important example of a log structure with an injective $u\:\calM\to\calO_X$ is as follows. Assume that $D\into X$ is a closed subset and set $\calM(\rmlog D)=\calO_{X_\et}\cap i_*\calO^\times_{U_\et}$ where $i\:U\into X$ is the open immersion of the complement $U=X\setminus D$. This is the log structure of elements invertible outside of $D$. Usually it is used when $D$ is the underlying closed set of a Cartier divisor; in this case $D$ is determined by the log structure and we call $\calM(\rmlog D)$ a {\em divisorial} log structure. In particular, the toroidal scheme structure $(X,U)$ can also be encoded in the log structure of $D$-monomial elements, where $D=X\setminus U$.

(3) The other extreme case is provided by so-called hollow log schemes with $u^m=0$ for any $m\in\calM\setminus\calM^\times$. Usually they show up when one restricts a log structure on $X$ onto a closed subscheme $Z$ such that $\calM_X$ is generically non-trivial on $Z$. Often this is the most economical way to encode certain information about the ambient scheme $X$ on $Z$. In particular, this is useful in deformation theory. In fact, it was such kind of an example, and not toroidal schemes, which led Fontaine and Illusie to introduce log schemes.
\end{exam}

\subsubsection{Associated log structure}
Any prelog structure $\calM$ can be canonically transformed into a log structure $\calM^a$. The idea is that $\calM$ is a log structure if and only if $u^{-1}(\calO_{X_\et}^\times)\toisom\calO_{X_\et}^\times$ so we should force this map to be an isomorphism. This works for any topology, so for shortness we only consider the \'etale one.

\begin{exer}
Given a prelog structure $\calM$ on $X$ let $\calM^a$ be the pushout of the diagram
$$\calM\hookleftarrow u^{-1}(\calO_{X_\et})^\times\to\calO_{X_\et}^\times.$$
Show that $\calM^a$ is a log structure, $\calM\to \calM^a$ is the universal homomorphism of $\calM$ to a log structure and the functor $\calM\mapsto \calM^a$ is left adjoint to the embedding of the category of log structures into the category of prelog structures.
\end{exer}

\begin{rem}
One can view the functor $\calM\mapsto \calM^a$ as an analogue of sheafification. Various operations on sheaves are defined by applying sheafification to the analogous operation in the category of presheaves (the sheafification is not needed for right exact functors, but is usually needed for other functors). In the same vein, various naive operations on log structures result in a prelog structure only (one works with sheaves, so the usual sheafification is used), and the functor $\calM\mapsto \calM^a$ should then be applied.
\end{rem}

\begin{defex}
Integralization $\calM^\int$ and saturation $\calM^\sat$ of a log structure $\calM$ are defined by applying the corresponding functors to monoids of sections, sheafification and then applying the functor $a$. Check that the resulting logarithmic structure is indeed saturated or integral.
\end{defex}

\subsubsection{Coherent log structures}
One reason to consider prelog structures was already discussed -- they form a simpler category and various operations on log structures often have prelog structures as intermediate results. Another reason is that prelog structures are used as charts (see \S\ref{chartsec} below), usually of finite type, for log structures, which are often very large because of the invertible part.

\begin{defin}
(i) A prelog structure is called {\em constant} if it is the sheafification of a homomorphism $P\to\Gamma(\calO_X)$ for a monoid $P$. For shortness we will not distinguish the monoid $P$ and its sheafification.

(ii) A log structure $\calM$ is called {\em quasi-coherent} if \'etale-locally there exists a constant prelog structure $P$ such that $\calM=P^a$. If, in addition, $P$ can be chosen finitely generated, $\calM$ is called {\em coherent}. We warn the reader that this notion is not related to coherence and quasi-coherence of $\calO_X$-modules.

(iii) A coherent and integral (resp. saturated) logarithmic structure is called {\em fine} (resp. {\em fs}).
\end{defin}

The following result can be found, for example, in \cite[Corollary~II.2.3.6]{logbook}

\begin{exer}
Show that a log structure is fine (resp. saturated) if and only if \'etale locally it possesses a fine (resp. fs) chart $P\to\calO_X$.
\end{exer}

Finally, for a morphism $f\:Y\to X$ one defines direct and inverse images of the log structures $\calM$ on $X$ and $\calN$ on $Y$:

\begin{defex}
(i) Show that $f_*(\calN)\times_{f_*(\calO_Y)}\calO_X$ is a log structure denoted (by a slight abuse of notation) $f_*(\calN)$.

(ii) The inverse image $f^*(\calM)$ is the log structure associated to the prelog structure $f^{-1}(\calM)\to f^{-1}(\calO_X)\to\calO_Y$. Show that, as expected, $f^*$ is left adjoint to $f_*$.
\end{defex}

\subsection{Logarithmic schemes}
Now we can introduce the category of log schemes.

\begin{defin}
(i) A log scheme $X$ is a tuple $(\uX,\calM_X,u_X)$, where $\uX$ is a scheme called the {\em underlying} scheme and $u_X\:\calM_X\to\calO_\uX$ is a log structure on $X$. Usually we will use looser notation when this cannot lead to a confusion, e.g. write $u$ instead of $u_X$ or just omit it, write $\calO_X$ instead of $\calO_\uX$ or even denote the underlying scheme by the same letter $X$.

(ii) A log scheme $X$ is called {\em quasi-coherent}, {\em coherent}, {\em integral}, {\em fine}, {\em fs}, etc., if the log structure $\calM_X$ is quasi-coherent, coherent, integral, fine, fs, etc.

(iii) A morphism of log schemes $f\:Y\to X$ consists of the underlying morphism $\uf\:\uY\to\uX$ of schemes and a homomorphism of log structures $f^*\calM_X\to\calM_Y$ compatible with the structure homomorphisms $u_X$ and $u_Y$.
\end{defin}

\begin{rem}
(i) We warn the reader that the word ``integral'' becomes slightly overused, because the notion of integral schemes means something completely different in the theory of schemes. So one should use it carefully, to avoid misunderstandings. Typically, if needed, one stresses that the {\em underlying scheme} is integral.

(ii) Already in \cite{Kato-log} Kato noticed that non-integral log schemes are too pathological and mainly restricted consideration to the category of fine log schemes. Although some aspects of a more general theory were developed quite systematically in \cite{logbook}, most of studies are done in the generality of fine log schemes, and this indeed seems to be the best choice. The second popular choice, which is often used once one works with log blowings up, is to work with the subcategory of fs log schemes. One benefit of working with these categories is that integral or saturated fiber products often reveal a nicer behaviour. In particular, log blowings up are log \'etale monomorphisms, see Exercise~\ref{logblow2} and this fact essentially uses integralization.
\end{rem}

\begin{rem}
(i) In case of a log structure with an injective $\calM_Y\into\calO_Y$ a morphism $f\:Y\to X$ is uniquely determined by the underlying morphism of schemes, and $\uf$ extends to $f$ if the log structure on $Y$ is ``larger'' than the image of the log structure on $X$. If $\calM_Y=\calM(\rmlog E)$ and $\calM_X=\calM(\rmlog D)$ are divisorial, then this happens if and only if $f^{-1}(D)\subseteq E$.

(ii) In particular, providing a toroidal scheme $(X,U)$ with the divisorial log structure $\calM_X=\calM(\rmlog(X\setminus U))$ we obtain a fully faithful embedding of the category of toroidal schemes into the category of log schemes. Naturally, such $(X,\calM_X)$ will be called a {\em toroidal log scheme}.

(iii) The other extreme is when the log schemes are hollow: $u_X=0$, $u_Y=0$. In this case, $\uf$ extends to $f$ via any homomorphism $f^*(\calM_X)\to\calM_Y$.
\end{rem}

\subsubsection{Strict morphisms}
A very important class of morphisms are those that ``minimally modify the log structure''. In the case of divisorial log structures this just means that $E=f^{-1}(D)$ and in general:

\begin{defin}
A morphism $f\:Y\to X$ of log schemes is called {\em strict} if $f^*(\calM_X)=\calM_Y$.
\end{defin}

In particular, for a log scheme $X$ any morphism $\uY\to\uX$ can be uniquely enhanced to a strict morphism of log schemes.

\begin{rem}
Any morphism $Y\to X$ canonically factors into a composition $Y\stackrel{g}\to Z\stackrel{h}\to X$ such that $h$ is strict and $\uY=\uZ$. It is natural to view $h$ as a scheme-like morphism and $g$ as a morphism which only increases the log structure. However, this factorization is not especially useful. The reason is that $g$ is not a ``monoidal-like'' morphism, see \S\ref{splitsec} below.
\end{rem}

\begin{exam}
(i) In the category of log schemes strict closed immersions play the role of usual closed immersions in the category of schemes. Often strict closed immersions have hollow sources. A typical example is when $X$ is a toroidal log scheme and $\uY\into\uX$ is a toroidal subscheme.

(ii) The simplest and most important case is when $X$ is the affine line with marked origin (i.e. $\uX=\Spec(k[t])$ and $\calM_X=(\bbN \log t)^a$, where we denote the generator of the monoid by $\log t$ to stress that it is mapped to $t$ by $u_X$) and $\uy=\Spec(k[t]/(t))$ is the origin. The induced log structure is $\calM_y=k^\times\oplus\bbN\log t$ and we call $(y,\calM_y)$ the {\em standard log point}. It is an analogue of points with non-reduced scheme structure in the usual algebraic geometry, in particular, we will later see that it is not log smooth and has non-trivial log differentials coming from the log direction $\log t$.

Note also that any other point $x=\Spec(k[t]/(t-a))$, $a\in k^\times$ has the trivial induced log structure since $\bbN\log t$ is mapped to $\calO_x^\times$ and hence the functor $\calM\mapsto\calM^a$ shrinks $k^\times\oplus\bbN\log t$ to $\calM_x=k^\times$ by sending $\log t$ to $a$.

(iii) Similarly, for each $n\ge 1$ the thick point $\uY_n=\Spec(k[t]/t^{n+1})$ can be provided with the log structure induced by $\bbN\log t$. It is neither injective, nor hollow.

(iv) One can also consider other log points, for example, for a toric monoid $P$ the log structure induced by $P\stackrel{0}{\to} k$ corresponds to the origin of the toric scheme $\Spec(k[P])$.
\end{exam}

\subsubsection{Log rings}
For a local work with log schemes it is often convenient to use the following logarithmic spectrum construction.

\begin{defin}
(i) A log ring $(A,P,u)$ consists of a ring $A$, a monoid $P$ and a homomorphism $u\:P\to(A,\cdot)$.

(ii) The logarithmic spectrum of a log ring, usually denoted just by $X=\Spec(P\stackrel u\to A)$ is the underlying scheme $\uX=\Spec(A)$ provided with the log structure induced by the prelog structure $u$.
\end{defin}

\begin{rem}
(i) In a sense this is a prelog ring, but the notion of a log ring does not make too much sense - even if $P=\Gamma(\calM_X)$, this condition will be lost after localizations.

(ii) A log scheme is Zariski if and only if it is a spectrum of log rings Zariski-locally. For general log schemes this is only true \'etale-locally, so in this aspect they are analogous to algebraic spaces.
\end{rem}

\subsubsection{Charts}\label{chartsec}
The notion of charts of log structures can be adopted to the category of log schemes once one replaces a monoid $P$ by the associated log scheme.

\begin{defex}\label{chartdef}
(i) For a monoid $P$ let $\bfA_P$ denote the log scheme whose underlying scheme is $\Spec(\bbZ[P])$ and the log structure is induced by the homomorphism $P\into\bbZ[P]$. If we work over a base scheme, for example, $S=\Spec(k)$, then we will use the notation $\bfA_{S,P}=S[P]:=S\times\bbZ[P]$.

(ii) Let $X$ be a log scheme. Show that giving a global chart $\phi\:P\to\calO_X$ for the log structure $\calM_X$ is equivalent to giving a strict morphism of log schemes $f\:X\to\bfA_P$. Any such morphism $f$ is called a global {\em chart} of $X$, and we will not make a real distinction between two presentations of a chart. A chart is called fine, saturated, sharp, etc. if the monoid $P$ is fine, saturated, sharp, etc. In particular, check that a coherent (resp. fine, resp. fs) log scheme is a log scheme which \'etale-locally possesses a finitely generated (resp. fine, resp. fs) chart.
\end{defex}

\begin{rem}
It is important to consider charts with non-sharp $P$ because various operations can produce non-sharp monoids. For example, removing the origin from $S[\bbN]=\Spec(k[t])$, where $S=\Spec(k)$, one obtains the scheme $S[\bbZ]=\Spec(k[t^{\pm 1}])$, which is also a chart of itself. On the other hand, its log structure is trivial, hence $S[\bbZ]\to S$ is also a chart.
\end{rem}

Nevertheless, when working locally at a point one might want to take a smallest possible chart and often this is possible. The follwoing notion is due to Kato:

\begin{defin}
A chart $P\to\calO_X$ is called {\em neat} at a geometric point $\ox\to X$ if $P\toisom\ocalM_\ox$. In particular, $P$ is automatically sharp.
\end{defin}

\begin{exam}
Let $X=\bfA_{k,\bbN}=\Spec(k[t])$. Then the tautological chart $\bbN\log t\to\calO_{X}$ is neat at the origin, but not at the other points, where a smaller chart (the trivial one) exists.
\end{exam}

The following example demonstrates one of rare aspects, in which fppf fine log schemes behave nicer than their \'etale analogues. However, even this is only needed when the log schemes are only fine but not fs. The following results can be found, for example, in \cite[\S II.2.3]{logbook}.

\begin{exer}\label{neatexer}
Let $X$ be a fine log scheme, let $\ox\to\uX$ be a geometric point and let $P=\ocalM_{X,\ox}$. Any neat chart at $\ox$ gives rise to a section $P\to\calM_{X,\ox}$ of the sharpening homomorphism.

(i) Show, that conversely any section $s\:P\to\calM_{X,\ox}$ of $\calM_{X,\ox}\to P$ induces a chart for a small enough \'etale neighborhood of $\ox$ and this chart is neat at $\ox$.

(ii) Furthermore, show that such a section exists if and only if the sequence $$1\to\calO^\times_{X_\et,\ox}\to\calM^\gp_{X,\ox}\to P^\gp\to 0$$ splits, and this is automatic whenever $P^\gp$ has no torsion of order divisible by $p=\cha(k(\ox))$ (\cite[Proposition~II.2.3.7]{logbook}).

(iii) Show that any fine fppf log scheme admits neat charts fppf-locally. (Hint: using the fppf topology one can also extract roots of order $p$.)

(iv) Show that if the log structure is fs, then $P^\gp$ is automatically torsion free, and hence even Zariski log schemes possess neat charts locally at points of $\uX$.
\end{exer}

\subsubsection{Monoidal morphisms}
We say that a morphism of log schemes $Y\to X$ is {\em monoidal} if \'etale-locally on $X$ it is the base change of morphisms of the form $\bfA_Q\to\bfA_P$. Informally speaking, $Y$ is obtained by first changing the monoidal structure and then adjusting the underlying scheme in the minimal needed way. Also, such morphisms can be viewed as base changes of morphisms of fans of monoids, e.g. see the informal notation \cite[Definition~9.10]{Kato-toric}. Main examples of such morphisms that are integralization, saturation, log blowings up and the morphism $\Log_X\to X$. They all be discussed later, and we start with the first two.

\begin{exer}\label{intsat}
(i) Show that for any coherent log scheme $X$ there exists a universal morphism $X\to X^\int$ (resp. $X\to X^\sat$) whose target is a fine (resp. fs) log scheme. In other words, the functor $X\mapsto X^\int$ (resp. $X\mapsto X^\sat$) is left adjoint to the embedding of the category of fine (resp. fs) integral schemes into the category of coherent log schemes. In addition, $X\to X^\int$ is a closed immersion and $X\to X^\sat$ is finite. (Hint: first, assuming that $X$ possesses a chart $X\to\Spec(\bbZ[P])$ show that $X^\int=X\otimes_{\bbZ[P]}\bbZ[P^\int]$ and $X^\sat=X\otimes_{\bbZ[P]}\bbZ[P^\sat]$ are as required. In general, use local-\'etale charts and \'etale descent of finiteness.)

(ii) Show that for any strict morphism of coherent schemes $Y\to X$ one has that $Y^\int=X^\int\times_XY$ and $Y^\sat=X^\sat\times_XY$.
\end{exer}


\subsubsection{Fiber products}
Similarly to the category of schemes, the category of log schemes and its subcategories possess all finite limits, of which we will use only fiber products.

\begin{exer}
(i) Let $\{X_i\}_{i\in I}$ be a finite diagram of log schemes. Show that $X=\lim_I X_i$ exists and can be described as follows: $\uX=\lim_I\uX_i$ and the log structure $\calM_X$ is the colimit of the pullbacks of $\calM_{X_i}$ to $\calM_X$. (In particular, this involves the functor $a$; the shortest way is to pullback $\calM_{X_i}$ as prelog structures, then take the colimit, and then apply $a$ once.)

(ii) Show that if all $X_i$ are fine or fs, then there exists a limit in the same category, and it coincides with $X^\int$ or $X^\sat$, respectively.
\end{exer}

\begin{rem}
Often one mentions $\int$ or $\sat$ in the notation of the limit, e.g. $Y\times_X^\sat Z$ or $(Y\times_XZ)^\int$, to avoid confusions. Sometimes, when only fs (resp. fine) log schemes are considered, this superscript can be omitted.
\end{rem}

To feel how this works, we start with the following almost tautological fact.

\begin{exer}
If $S$ is a scheme, $P\to Q$, $P\to R$ are homomorphisms of finitely generated monoids and $X=S[P]$, $Y=S[Q]$, $Z=S[R]$, then $Y\times_XZ=S[Q\oplus_P R]$, $(Y\times_XZ)^\int=S[(Q\oplus_PR)^\int]$ and $(Y\times_XZ)^\sat=S[(Q\oplus_PR)^\sat]$.
\end{exer}

Now we can consider an archetypical example of a subtle behaviour of fine log schemes, as opposed to coherent log schemes or schemes.

\begin{exer}\label{fineexer}
(i) Let $X=\Spec(k[x,y])$ with the log structure induced by $\bbN\rmlog(x)\oplus\bbN\rmlog(y)$, and let $Y=\Spec(k[x,\frac{y}{x}])$ with the log structure induced by $\bbN\rmlog (x)\oplus\bbN(\rmlog(y)-\rmlog(x))$. Show that $Y=Y\times^\int_XY$, in particular, the integralization functor cuts off the $\bbA^2$ component from $\uY\times_\uX\uY$.

(ii) More generally, assume that $P\subsetneq Q$ are toric monoids such that $P^\gp=Q^\gp$. Then $\bfA_Q=\bfA_Q\times^\int_{\bfA_P}\bfA_Q$, that is, the morphism $\bfA_Q\to\bfA_P$ is a monomorphism in the category of fine log schemes, but not in the category of coherent log schemes.
\end{exer}

\begin{rem}\label{finerem}
Usual blowing up have some nasty properties. For concreteness, consider the blowup up chart $f\:\uY=\Spec(k[x,y'=\frac yx])\to\uX=\Spec(k[x,y])$ at the origin. Clearly, $f$ is non-flat, the dimension of the fiber over the origin jumps, and this even gives rise to the new irreducible component $\Spec(k[y'_1,y'_2])$ in $\uY\times_\uX\uY$. Assume now that $\uX$ and $\uY$ are provided with the log structures generated by $\log(x),\log(y)$ and $\log(x),\log(y')$. Then the $\bbA^2$ component in the product acquires a non-cancellative log structure with monoid presented by generators $\log(x),\log(y'_1),\log(y'_2)$ and relation $\log(x)+\log(y'_1)=\log(y)=\log(x)+\log(y'_2)$. This forces the integralization functor to remove the $\bbA^2$ component and only its diagonal, which is the intersection with the other component, is left. In fact, we will later see that many similar morphisms (log blowings up and their charts) are monomorphisms in the fine category.
\end{rem}

And here is an archetypical example of a property of the saturated category. It explains why one usually restricts the setting to fs log schemes when studying Kummer covers.

\begin{exer}\label{kumexer}
(i) Let $X=\Spec(k[x])$ with the log structure $\bbN\log(x)$ and $Y=\Spec(k[x^{1/2}])$ with the log structure $\frac{1}{2}\bbN\log(x)$. We will later see that the Kummer cover $Y\to X$ is log \'etale when $\cha(k)\neq 2$. Check that $Z=Y\times_XY=(Y\times_XY)^\int$ contains two components (diagonal and antidiagonal) intersecting over the origin and the characteristic $\ocalM_z$ at the intersection point $z$ is the non saturated monoid $\frac{1}{2}\bbN\oplus_\bbN\frac{1}{2}\bbN$ with generators $(\frac 12,0),(0,\frac 12)$ subject to the relation $(1,0)=(0,1)$. It is isomorphic to the submonoid of $\bbZ/2\bbZ\oplus\bbN$ obtained by removing the element $(1,0)$.

(ii) Show that $Z^\sat$ is the normalization of $Z$ and it is just the disjoit union of the two copies of $Y$ -- the diagonal and the antidiagonal. The characteristic of the two points over $z$ is $\bbN$ -- the sharpening of $\ocalM_z^\sat\toisom\bbZ/2\bbZ\oplus\bbN$.

(iii) More generally, if $S$ is a scheme, $P\subseteq Q$ is a Kummer extension of toric monoids, i.e. $Q$ is the saturation of $P$ in $Q^\gp$, and $X=S[P], Y=S[Q]$, then $Z=Y\times_X^\sat Y$ is isomorphic to the split cover $Y\times G$ for $G=Q^\gp/P^\gp$. So, the log \'etale $G$-Galois cover $S[Q]\to S[P]$ behaves similarly to \'etale covers only in the category of fs log schemes.
\end{exer}

\subsection{Logarithmic regularity}
In this section we only consider fs log schemes.

\subsubsection{The logarithmic stratification}
Each log scheme $X$ possesses a natural stratification by $\rk(\ocalM_x)=\dim_\bbQ(\ocalM_x^\gp\otimes\bbQ)$.

\begin{exer}
(i) Let $C_{\ge d}\subseteq\uX$ be the locus on which the rank of the characteristic monoid is at most $d$. Show that these sets are closed and hence induce a stratification by the locally closed sets $C_d=C_{\ge d}\setminus C_{>d}$.

(ii) Refine this stratification to a {\em log stratification} $\uX_{d}$ of $X$ by (non-necessarily reduced) locally closed subschemes as follows: if the structure is Zariski at $x$ and $d=\rk(\ocalM_x)$, then the stratum $\uX_{\ge d}$ at $x$ is given by the vanishing of $u^{\calM_x^+}$, where $\calM_x^+$ is the maximal ideal of $\calM_x$. Show that this is compatible with strict morphisms and hence extends to arbitrary log structures by \'etale descent. Finally, set $X_{d}=X_{\ge d}\setminus X_{>d}$

(iii) Show that, indeed, $C_{\ge d}$ is the reduction of $X_{\ge d}$. Also, show that the log strata of the chart log schemes $\bfA_P$ are reduced.
\end{exer}

\subsubsection{Logarithmic regularity}
The following definition is a far-reaching generalization of the classical fact recalled in Exercise~\ref{sncexer}(iii).

\begin{defin}
A locally noetherian fs logarithmic scheme $X$ is {\em logarithmically regular} if each locally closed subscheme $\uX_{r}=\uX_{\ge r}\setminus\uX_{\ge r+1}$ is regular (in particular, reduced) and of codimension $r$.
\end{defin}

\begin{rem}
The original definition by Kato only considered fs log schemes. Most of results about log regularity can be extended to fine log schemes, and this was worked out by Gabber, but we will not touch this direction in the notes.
\end{rem}

\subsubsection{Log parameters}
As in the case of regular schemes, when working with log regular log schemes it is very convenient to use local parameters. The classical notion is generalized as follows:

\begin{defin}
Let $X$ be a log regular log scheme with a geometric point $\ox\to X$ over $x\in X$. Let $r=\rk(\ocalM_\ox)$. By a {\em regular family of parameters} at $\ox$ we mean a section $s\:\ocalM_\ox\to\calM_\ox$ and elements $t_1\..t_d\in\calO_x$, such that $t_1\..t_d$ restrict to a regular family of parameters of $\uX_r$ at $x$. In particular, $d$ is the dimension of $\uX_r$ at $x$ and $r+d=\dim_x(X)$. We call $t_i$ {\em regular parameters} and the elements $u^{s(m)}$ for $0\neq m\in\ocalM_x$ will be called {\em logarithmic parameters}.
\end{defin}

In view of Exercise \ref{neatexer}(iv) such families of parameters exist, and if the log structure at $x$ is Zariski, one can even construct $s$ Zariski locally. As in the classical case, it follows from the definition that any regular family of parameters generates the maximal ideal at $\ox$. Furthermore, parameters naturally give rise to very explicit \'etale and formal charts.

\begin{exer}\label{logregex}
Assume that a log scheme $X$ is Zariski and log regular at a point $x\in X$ with $P=\ocalM_x$, and let $s\:P\to\calM_x$ and $t_1\..t_d\in\calO_x$ be a regular family of parameters at $x$. If $k=k(x)$ is of positive characteristic $p$, let $C_k$ be a Cohen ring of $k$, that is a DVR with maximal ideal $(p)$ and residue field $k$. By Cohen's theorem if $X$ is of equal characteristic at $p$, then there exists a field of coefficients $i\:k\into\hatcalO_x$, while in the mixed characteristic case there exists a ring of coefficients $i\:C_k\into\hatcalO_x$. Prove the following theorem of Kato (see \cite[Theorem~3.2]{Kato-toric}), where $A\llbracket P\rrbracket$ denotes the formal completion of $A[P]$ at the ideal $A[P^+]$:

(i) In the equal characteristic case the natural homomorphism $k\llbracket P\rrbracket\llbracket t_1\..t_d\rrbracket\to\hatcalO_x$, induced by $i$ and the parameters, is an isomorphism.

(ii) If $\cha(k)>0$, then the natural homomorphism $C_k\llbracket P\rrbracket\llbracket t_1\..t_d\rrbracket\to\hatcalO_x$ induced by $i$ and the parameters is surjective with a principal kernel $(\theta)$, where $\theta\equiv p$ modulo $(P^+,t_1\..t_d,p)^2$.
\end{exer}

\begin{rem}
(i) The equal characteristic $p$ case naturally shows up in both cases, since one can take $\theta=p$.

(ii) In the equal characteristic case Kato's theorem tells that log regularity is the same as being formally-locally isomorphic to a toric variety. So, the theory of log regular schemes can be viewed as a generalization of toroidal geometry to the mixed characteristic case.

(iii) A relatively difficult theorem asserts that log regularity is preserved by localizations. Kato's original proof is incomplete, but Gabber later provided missing arguments. The source of the difficulty is clear -- regular parameters at generizations of $x$ are not related to parameters at $x$. In fact, this is completely parallel to the situation with usual regularity. However, in the classical case there is a conceptual proof which uses Serre's cohomological criterion of regularity, and no logarithmic analogue was found so far.
\end{rem}

\subsubsection{Log regularity and toroidal varieties}
Finally, let us describe log regular log varieties over a field $k$. As in the case of usual varieties, a nice description (in terms of smoothness) is possible only for simple points $x$, that is, points for which the extension $k(x)/k$ is separable.

\begin{exer}
Assume that $X$ is a log variety over $S=\Spec(k)$ and $x\in X$ is such that $k(x)/k$ is separable and $X$ is log regular at $x$. Let $\ox\to X$ be a geometric point over $x$.

(i) Prove that locally at $\ox$ there exist a chart $f\:U\to S[P\oplus\bbZ^d]$ with an \'etale $f$, where $P=\ocalM_\ox$ and $d$ is the dimension of the logarithmic stratum at $x$. Deduce that $X$ is toroidal at $\ox$, and thus being log regular and toroidal at a simple point are equivalent. (Hint: use parameters and work with \'etale topology instead of the formal one.)

(ii) Prove that any chart $f\:U\to S[P]$, which is neat at $\ox$, is smooth at the image of $\ox$.
\end{exer}

\section{Morphisms of logarithmic schemes}\label{morsec}
Our next goal is to study morphisms of log schemes in more detail.

\subsection{Charts}
Recall that charts of log schemes were defined in Definition \ref{chartdef}. Naturally, by a chart of a morphism of log schemes one means compatible charts of the source and the target:

\begin{defin}
Let $f\:Y\to X$ be a morphism of log schemes. A {\em chart} of $f$ consists of charts $Y\to\bfA_Q$ and $X\to\bfA_P$ and a homomorphism $\phi\:P\to Q$ such that the compositions $P\to\Gamma(\calM_X)\to\Gamma(\calM_Y)$ and $P\to Q\to\Gamma(\calM_Y)$ coincide. Equivalently, the chart is a commutative diagram
$$
\xymatrix{
Y\ar[d]^f\ar[r] & \bfA_Q\ar[d]^{\bfA_\phi}\\
X\ar[r] & \bfA_P
}
$$
whose horizontal lines are charts. One says that the chart is {\em modeled on} $\phi$.
\end{defin}

It is easy to construct charts, morally, one just starts with a chart $P\to X$ and enlarges it by adding enough elements of $\calM_Y$.

\begin{exer}
Let $f\:Y\to X$ be a morphism of fine log schemes, $\ox\to X$ a geometric point and $\oy\to Y$ a geometric point above $\ox$. Prove that any \'etale-local fine chart $P\to\calM_\ox$ at $\ox$ extends to a chart of $f$ \'etale-locally at $\oy$. (Hint: start with any fine chart $Q'\to\calM_\oy$ and take $Q$ to be the image of $P\oplus Q'\to\calM_\oy$.)
\end{exer}

\subsubsection{The standard splitting}\label{splitsec}
If a log scheme $X$ is provided with a chart $X\to\bfA_P$ and $P\to Q$ is a homomorphism of monoids we will use the notation $X_P[Q]=X\times_{\bfA_P}\bfA_Q$, which indicates that $X_P[Q]$ is obtained by a ``base change of the monoidal structure''.
A chart of $f$ induces a very useful splitting of $f$ into the composition $Y\to X_P[Q]\to X$, where the first morphism is strict, as $\bfA_Q$ is a chart of both the source and the target, and the second morphism is monoidal. In a sense, this decomposes $f$ into a composition of a scheme-like morphism and a monoidal-like morphism, and such a splitting is ubiquitous in log geometry. However, it is non-canonical and exists only \'etale-locally.

\subsubsection{Neat charts}
As with the neat absolute charts, to efficiently work with charts of morphisms one would like to construct minimal charts, or even a minimal chart extending a given neat chart of the target. It turns out that in general one cannot construct charts modeled on homomorphisms of sharp monoids, as we are going to demonstrate. Let $f\:Y\to X$ be a morphism of Zariski log schemes, $y\in Y$, $x=f(y)$, $Q=\ocalM_y$ and $P=\ocalM_x$. The induced homomorphism $\phi\:P\to Q$ has no kernel, since any non-unit of $\calM_x$ is mapped to a non-unit in $\calO_x$ and hence also to a non-unit in $\calO_y$. Nevertheless, it might happen that $\phi$ is not injective, but $\Ker(\phi^\gp)\cap P=0$. Here is an archetypical example:

\begin{exam}
Let $X=\Spec(k[s,t])$ with the log structure generated by $s$ and $t$, and $Y=\Spec(k[s,u])$ with $u=t/s$ a chart of the blowing up of $X$ at the origin $x\in X$ with the log structure generated by $s$ and $u$. Let $y\in Y$ be a point given by $s=0$, $u=a\in k^\times$. Then $P=\ocalM_x=\bbN\log(s)\oplus\bbN\log(t)$ and $Q=\ocalM_y=\bbN\log(s)$ because $u\in\calO_y^\times$. The map $P\to Q$ sends $\log(t)$ to $\log(s)$ and $\log(u)=\log(t)-\log(s)$ generates the kernel of the map $P^\gp\to Q^\gp$.
\end{exam}

The above example motivates the following definition with the idea to provide $Q$ with the minimal amount of units so that the homomorphism $P\to Q$ is injective.

\begin{defin}
(i) The {\em relative characteristic monoid} of $f\:Y\to X$ is $\calM_{Y/X}=\Coker(f^*(\calM_X)\to\calM_Y)$. Note that it also coincides with $\Coker(f^{-1}(\ocalM_X)\to\ocalM_Y)$, hence $\calM_{Y/X}=\ocalM_{Y/X}$.

(ii) Assume that $Y\to\bfA_Q$, $X\to\bfA_P$, $\bfA_\phi$ is a chart of a morphism of log schemes $f\:Y\to X$, and $y\in Y$ is a point with $x=f(y)$. The chart is called {\em neat} at $y$ if $\phi$ is injective and the induced homomorphism $\Coker(\phi^\gp)\to\calM^\gp_{Y/X,y}$ is an isomorphism.
\end{defin}

Neat charts always exist in fppf topology, while in \'etale topology there might be an obstacle if $\Coker(\phi^\gp)$ has a $p$-torsion for $p=\cha(k(y))$. Proof of the following result, which can be found in \cite[Theorem~II.2.4.4]{logbook}, involves a bit of diagram chasing with monoids and their Grothendieck groups.

\begin{exer}
Assume that $f\:Y\to X$ is a morphism of log schemes, $\oy\to Y$ a geometric point whose image in $X$ is $x\in X$. Assume that $h\:U\to\bfA_P$ is a neat chart of a neighborhood of $x$. If $\Ext^1(\calM^\gp_{Y/X,\oy},\calO_\oy^\times)=0$, then locally along $\oy$ there exists a neat chart of $f$ that extends $h$. In particular, a neat chart exists when the order of torsion of $\Coker(\calM^\gp_{Y/X,\oy})$ is invertible in $k(\oy)$.
\end{exer}

\begin{rem}
The same argument proves that neat charts always exist fppf locally. Similarly, if the log structure at $y$ is Zariski and $\calM^\gp_{Y/X,\oy}$ is torsion free, a neat chart exists Zariski locally, but if the relative characteristic monoid contains a torsion, one might need to extract roots of some units, ending up with an \'etale-local or fppf-local chart -- depending on the torsion.
\end{rem}

\subsection{Logarithmic smoothness}

\subsubsection{Log thickenings}
The following definition is a direct extension of its scheme-theoretic analogue.

\begin{defin}
(i) A {\em log thickening} is a strict closed immersion $S\into T$ given by a nilideal $\calI\subset\calO_T$.

(ii) A morphism of log schemes $f\:Y\to X$ is called {\em formally log smooth} (resp. {\em formally log \'etale}, resp. {\em formally log unramified}) if for any log thickening and compatible morphisms $i\:S\to Y$, $T\to X$ \'etale locally on $T$ there exists (resp. there exists unique, resp. there exists at most one) lifting $T\to X$ making the diagram commutative:
$$
\xymatrix{
S\ar@{^{(}->}[d]^i\ar[r] & Y\ar[d]^f\\
T\ar[r]\ar@{.>}[ru] & X
}
$$
(iii) A morphism $f$ is {\em log unramified} if it is formally log unramified and $\uf$ is of finite type. A morphism $f$ is {\em log smooth} (resp. {\em log \'etale}) if it is formally log smooth (resp. log \'etale) and $\uf$ is finitely presented.
\end{defin}

Similarly to the theory of schemes one can study these notions most effectively by use of log differentials.

\subsubsection{Logarithmic derivations}
Log geometry has a version of the theory of derivations and differentials. The idea is to extend the usual theory by adding logarithmic differentials of monomials.

\begin{defin}
Let $A\to B$ be a homomorphism of log rings $A=(P\to\uA)$ and $B=(Q\to\uB)$ and let $N$ be a $\uB$-module. An {\em $A$-log derivation} $(d,\delta)\:B\to N$ consists of an $\uA$-derivation $\uB\to N$ and a homomorphism $\delta\:Q\to N$ such that $\delta(P)=0$ and $\delta(q)=\frac{d(u^q)}{u^q}$. The $\uB$-module of all log derivations will be denoted $\Der_{B/A}(E)$.
\end{defin}

\begin{rem}
One should view $\delta$ as derivation of the branch of the logarithm of the monomial $u^q$ specified by $q$.
\end{rem}

\subsubsection{Logarithmic differentials}
Similarly to usual K\"ahler differentials, there always exists a universal log derivation whose target is the module $\Omega^1_{B/A}$ of {\em logarithmic differentials}. This is not so standard, but we omit 1 in the notation because we will never consider modules of differential $p$-forms with $p>1$.

\begin{exer}
Prove that the functor $\Der_{B/A}(\cdot)$ is representable and describe the corresponding universal module as follows: $\Omega^1_{B/A}$ is the quotient of $\Omega^1_{\uB/\uA}\oplus(B\otimes(Q^\gp/P^\gp))$ by the submodule generated by the relations $(du^q,-u^q\otimes q)$ with $q\in Q$, where we denote the image in $Q^\gp/P^\gp$ also by $q$.
\end{exer}

We will not go into details, but most of classical results about K\"ahler differentials, such as base change, compatibility with localizations and fundamental sequences, extend to log rings and log differentials. In addition, the module of differentials of any strict \'etale homomorphism vanishes, hence by use of \'etale descent the definition can be globalized to the definition of an $\calO_Y$-module $\Omega^1_{Y/X}$ for any morphism $Y\to X$ of log schemes, and it comes equipped with the universal log derivation $(d,\delta)\:(\calM_X\to\calO_X)\to\Omega^1_{Y/X}$. In particular, $d$ induces a homomorphism $\Omega^1_{\uY/\uX}\to\Omega^1_{Y/X}$. As usual in log geometry, to get some feeling of a new notion one should look at the two polar cases, and the case of strict morphisms reduces to the usual scheme theory.

\begin{exex}
Use the previous exercise to prove that:

(i) If $Y\to X$ is a strict morphism of log schemes, then $\Omega^1_{\uY/\uX}=\Omega^1_{Y/X}$.

(ii) If $R$ is a ring, $Y=\bfA_{R,P}$, $Y=\bfA_{R,Q}$ and $f=\bfA_{R,\phi}$ for $\phi\:P\to Q$, then $$\Omega^1_{Y/X}=\Coker(\phi^\gp)\otimes R[Q].$$
In particular, if $\bbQ\subseteq R$, then $\Omega^1_{Y/X}$ is free and its rank is the rank of $\Coker(\phi^\gp)$, and if $\bfF_p\subseteq R$, then $\Omega^1_{Y/X}$ is free and its rank is the $p$-rank of $\Coker(\phi^\gp)$. (Hint: show that $\Der_{R[Q]/R[P]}(E)=\Hom(\Coker(\phi^\gp),E)$.)

(iii) By a slight abuse of notation for a sharp monoid $P$ and a ring $R$ we denote by $P\stackrel 0\to R$ the homomorphism of monoids taking $P^+$ to 0 (but 0 goes to 1). Compute the differentials of a log point $\Spec(P\stackrel 0\to k)$. More generally, show that for any ring $R$ and the hollow log scheme  $X=\Spec(P\stackrel 0\to R)$ there is an isomorphism of $R$-modules $\Omega^1_{X/R}=P^\gp\otimes R$.
\end{exex}

\subsubsection{Chart criterion}
The main theorem of Kato about logarithmic smoothness gives the following criterion for log smoothness in terms of charts, see \cite[Theorem~3.5]{Kato-log}:

\begin{theor}\label{smoothchart}
For a morphism $f\:Y\to X$ of fine log schemes and a geometric point $\oy\to\uY$ the following conditions are equivalent:
\begin{itemize}
\item[(i)] $f$ is log smooth (reps. log \'etale) at $\oy$,

\item[(ii)] Locally at $\oy$ there exists a chart $V\to\bfA_Q$, $U\to\bfA_P$ modelled on $\phi\:P\to Q$ such that
\begin{itemize}
\item[(a)] The morphism $V\to U_P[Q]$ is smooth (resp. \'etale)

\item[(b)] $\Ker(\phi^\gp)$ and $\Coker(\phi^\gp)_\tor$ (resp. $\Ker(\phi^\gp)$ and $\Coker(\phi^\gp)$) are finite of order invertible in $k(\oy)$.
\end{itemize}

\item[(iii)] The same condition as (ii) but with $V\to U_P[Q]$ \'etale.
\end{itemize}
\end{theor}

\begin{exer}
(i) Check that the condition on $\phi$ is equivalent to smoothness (resp. \'etaleness) of the morphism of diagonalizable groups $D_{k(\oy),Q^\gp}\to D_{k(\oy),P^\gp}$ and by the above criteria it is also equivalent to smoothness (resp. \'etaleness) of the map $\bfA_{k,Q}\to\bfA_{k,P}$.

(ii) Deduce (iii) from (ii) by enlarging the chart to $Q\oplus\bbZ^n$, where $t_1\..t_n$ are regular parameters on the fiber of $V\to U_P[Q]$ through $\oy$. (Hint: for example, one can send the generators of $\bbZ^n$ to $1+t_i$ (i.e. one increases $Q$ by adding the elements $\log(1+t_i)$.)

(iii) Check that in condition (b) one can also achieve that $\Ker(\phi)=0$ (Hint: again, just increase $Q$ accordingly.)
\end{exer}

Here as a very typical example involving log points. We will return to this setting also in example~\ref{blowexam}.

\begin{exam}
(i) Show that log points $\Spec(P\stackrel 0\to k)$ with a sharp fine $P\neq 0$ are not log smooth over $k$.

(ii) Take $P=\bbN\log(t)$ and consider the nodal curve $C=\Spec(k[x,y])$ with the log structure induced by $Q=\bbN\log(x)\oplus\bbN\log(y)$. Show that the morphism $C\to P$ taking $\log(t)$ to $a\log(x)+b\log(y)$ with $a,b\in\bbN$ is log smooth if and only if $(a,b)$ is invertible in $k$. Show that if $(a,b)\in k^\times$, then (unlike $\Omega^1_{\underline{C}/k}$) the module $\Omega^1_{C/P}$ is invertible; in fact, it is generated by $\delta(x),\delta(y)$ subject to the relation $a\delta(x)+b\delta(y)=0$.

(iii) Give another proof of the results of (ii) by noting that the morphism $C\to P$ is just the fiber over the origin of the morphism $\bfA_Q\to\bfA_P$ given by $t=x^ay^b$.
\end{exam}

An important theorem of Kato states that log regularity is preserved by log smooth morphisms, see \cite[Theorem~8.2]{Kato-toric}:

\begin{theor}\label{regsmooth}
If $X$ is a log regular log scheme and $Y\to X$ is log smooth, then $Y$ is log regular.
\end{theor}

\subsection{Logarithmic blowings up}

\subsubsection{Log ideals}

\begin{defin}
Let $X$ be a fine log scheme. By a {\em log ideal} we mean a coherent sheaf of ideals $J\subseteq\calM_X$, where coherence means that locally around any geometric point $\ox$ the ideal is generated by $J_\ox$. A log ideal is called {\em invertible} if it is locally generated by a single element.
\end{defin}

\begin{rem}
Since monoids $\calM_X(U)$ are integral there is no need to impose a condition analogous to being a non-zero divisor in a ring, and invertible ideals are preserved by arbitrary pullbacks, unlike the theory of schemes.
\end{rem}

\subsubsection{Log blowings up}
Log blowings up are defined by the same universal property as usual blowings up, but with log ideals used instead of ideals.

\begin{defin}
Let $X$ be a log scheme and $J$ a log ideal on $X$. A morphism $f\:Y\to X$ is called the {\em log blowing up} of $X$ along $J$ and denoted $LogBl_J(X)\to X$ if $f$ is the universal morphism of fine log schemes such that $f^{-1}(J)$ is invertible. The saturated log blowing up is defined by the same property but in the category of fs log schemes. Thus, it is nothing else but the saturation of the log blowing up.
\end{defin}

\begin{exer}\label{logblow1}
Log blowings up are preserved by any base changes $h\:X'\to X$, i.e. $Y'=LogBl_{h^{-1}J}(X')=Y\times_XX'$. In particular, applying this to $Y$ one obtains that $Y\times_XY=Y$ and hence $Y\to X$ is a monomorphism.
\end{exer}

Now we are going to prove that log blowings up exist and are, in fact, monoidal morphisms, as one might expect, since they realize a monoidal construction.

\begin{exer}\label{logblow2}
We will construct $Y=LogBl_J(X)$ in a few steps.

(i) Prove that if $X=\bfA_P$ and $I=\bbZ[J]$, then $\uY=\Bl_I(X)$ and the log structure on the $s$-chart $\uY_s=\Spec(\bbZ[P][\frac{I}{s}])$ of the blowing up is given by $P[J-s]$ (the submonoid of $P^\gp$ generated by $P$ and the elements of $J-s$). Furthermore, the chart $Y_s$ of $Y$ is the universal log scheme over $X$ such that the pullback of $J$ to $Y_s$ is principal generated by $s$.

(ii) Prove that if $X$ possesses a chart $X\to\bfA_P$ and $J$ is generated by an ideal $J_0\subseteq P$, then $LogBl_J(X)=LogBl_{J_0}(\bfA_P)\times_{\bfA_P}X$. Deduce that log blowing up $Y\to X$ always exists and is a proper monoidal morphism. (Hint: use part (i), Exercise~\ref{logblow1} for the first claim and then apply \'etale descent.)

(iii) Prove that any log blowing up is a log \'etale morphism. (Hint: by base change and \'etale descent this reduces to the case described in (i) and then Kato's criterion does the job.)

(iv) Finally, reinterpret examples we have seen: using Exercise~\ref{fineexer}(ii) prove directly that $Y=Y\times_XY$, and explain which chart of  a log blowing up is considered in Remark~\ref{finerem}.
\end{exer}

Typically, one asks when a log blowing up coincides with the usual blowing up along the ideal $u^J\calO_X$ generated by $J$, because this seems to be the most adequate situation. Clearly, this happens if and only if the relevant toroidal blowings up of the charts $\bfA_P$ are compatible with the chart maps $h\:U\to\bfA_P$. In particular, this is automatic if $h$ is flat or, at least, what is called $\Tor$-independent from $J$, that is, $\Tor_1^{\bbZ[P]}(\calO_X,\bbZ[P]/J^n\bbZ[P]])=0$. Using the latter criterion Nizio{\l} proved the following claim, see \cite[Proposition~4.3]{Niziol}:

\begin{theor}\label{logblowth}
If $X$ is a log regular log scheme and $J$ is a log ideal, then the underlying scheme of the $LogBl_J(X)$ is the blowing up of $X$ along the induced ideal $u^J\calO_X$. In addition, the saturated log blowing up $LogBl_J(X)^\sat$ is log regular and the log structure is induced by the union of the preimage of the toroidal divisor on $X$ and the exceptional divisor.
\end{theor}

In fact the approach with Tor functors is only needed in the more difficult case of mixed characteristics.

\begin{exer}
Prove Theorem \ref{logblowth} when $X$ is equicharacteristic. (Hint: formal completions of noetherian schemes are flat and hence compatible with blowings up. Formally locally $X$ looks as $\Spec(k\llbracket P\rrbracket\llbracket t_1\.. t_d\rrbracket)$, hence the blowing up is described easily via the base change from the toric case via the flat homomorphism $k[P]\into k\llbracket P\rrbracket\llbracket t_1\.. t_d\rrbracket$. The log regularity follows from Theorem~\ref{regsmooth}.)
\end{exer}

Now let us discuss a few cases, that are simpler to compute but are often viewed as pathological and not worth consideration. In particular, one can easily have that the proper morphism $LogBl_J(X)\to X$ is not birational. Using such morphisms becomes critical if one wants to study relative desingularization over singular bases (e.g. the log point or thick log point with a non-reduced scheme structure), but this direction has not studied yet in the literature, and we will only mention it in a couple of remarks later.

\begin{exex}\label{blowexam}
(i) Show that the log blowing up of the ideal $P^+$ on the log point $\Spec(P\stackrel 0\to k)$ modeled on a fine sharp monoid $P$ is of dimension $\rk(P)-1$.

(ii) Let $s=\Spec(\bbN\log(t)\stackrel 0\to k)$ be the standard log point (i.e. $t=0$ in $k$) and $$C_s=\Spec(\bbN\log(x)\oplus\bbN\log(y)\to k[x,y]/(xy))$$ a log smooth $s$-curve with $\log(t)$ mapped to $\log(x)+\log(y)$. Let $J$ be the maximal ideal of $\bbN\log(x)\oplus\bbN\log(y)$. Show by a direct computation that $X_s=LogBl_J(C_s)$ is also log smooth over $s$, the map $X_s\to C_s$ is an isomorphism over the complement of the origin $O\in C_s$ and the preimage $E$ of $O$ is a non-reduced double $\bfP^1_k$ component with the nilpotent ideal $(\veps)$ and the log structure given by $\bbN\log(\veps)$ and $\log(t)$ mapping to $2\log(z)$.

(iii) Now, embed $s$ as the origin of $S=\Spec(\bbN\log(t)\to k[t])$ and $C_s$ as the closed fiber of a log smooth (even semistable) $S$-curve $$C=\Spec(\bbN\log(x)\oplus\bbN\log(y)\to k[x,y])$$ with $t=xy$. Recall that $LogBl_J(C_s)$ is the pullback of the log scheme $LogBl_J(C)$, whose underlying scheme is just $Bl_O(C)$. Use this to conceptually explain the results of (ii), in particular, the reason why the new component is doubled (has a non-reduced structure).
\end{exex}

\subsection{Logarithmic \'etaleness}
In this section we restrict to the fs setting, in which Kummer covers are usually studied. We just give definitions, check simplest properties and mention various directions studied in the literature.

\subsubsection{Kummer \'etale morphisms}
A homomorphism of toric monoids $\phi\:P\to Q$ is called {\em Kummer} if $P^\gp\subseteq Q^\gp$ is of finite index and $Q$ is the saturation of $P$ in $Q^\gp$. A log \'etale morphism of log schemes $Y\to X$ is called {\em Kummer} if the induced homomorphisms of monoids $\ocalM_x\to\ocalM_y$ are Kummer. A Kummer \'etale cover $Y\to X$ is a surjective Kummer \'etale morphism.

\begin{exer}
Check that, indeed, this notion of a covering defines a Grothendieck topology called Kummer \'etale of $\ket$ topology. (Hint: this mainly reduces to the check that Kummer \'etale covers are preserved by base changes.)
\end{exer}

\begin{rem}
(i) Kummer \'etale topology is the closest analogue of the \'etale topology in the setting of fs log schemes. For example, see Exercise~\ref{kumexer}. An important fact is that the theory of descent works pretty similarly to the case of \'etale (or flat) topology: $\calO_{X_\ket}$ and, more generally, representable functors are sheaves (see \cite[Proposition~2.18 and Theorem~2.20]{Niziol-K-theory-of-log-schemes-I}), etc. Ideals in $\calO_{X_\ket}$ are called {\em Kummer ideals}. Working with them provides a convenient formalism for extracting roots from monomials.

(ii) Theories of Kummer \'etale and general log \'etale cohomologies are now developed rather deeply, see \cite{Nakayama1} and \cite{Nakayama2}. The first one is simpler, but in order to have some fundamental theorems in full generality one has to work with the whole log \'etale site.
\end{rem}

\subsubsection{Log \'etale site}
For the sake of completeness, let us discuss how one defines the notion of a log \'etale covering in general. Our motivation is just to see a few more examples from log geometry, and we will not discuss the log \'etale cohomology theories.

\begin{exer}
Construct an example of surjective log \'etale morphisms $Y\to X$ and $Z\to X$ such that $Y\times_XZ\to X$ is not surjective. In particular, a base change does not have to be surjective, and hence can even be empty. (Hint: for example, one can take the plane with the monoid $\bbN^2$, apply log blowing up to the origin and another log blowing up to one of the two preimages of the origin with characteristic $\bbN^2$, obtaining a log \'etale morphism $X'\to X$ with the exceptional divisor consisting of two components $E=E_1\cup E_2$. Then $Y=X'\setminus E_1$ and $Z=X'\setminus E_2$ do the job. It is also instructive to consider a purely combinatorial (or toric) description of this example.)
\end{exer}

The above example shows that one should be careful with the notion of surjectivity. Naturally, we would like to declare any log blowing up $Y\to X$ to be a cover, but for any $Y'\subsetneq Y$ the morphism $Y'\to X$ should not be a cover. So, one defines the log \'etale topology to be the topology generated by Kummer \'etale covers and log blowings up.

\begin{exer}
Let $f\:Y\to X$ be a log \'etale morphism. Show that $f$ is a log \'etale cover if and only if for any log blowing up $X'\to X$ the morphism $Y\times_XX'\to X'$ is surjective.
\end{exer}

\section{The stacks $Log_S$}\label{Logsec}
This section is devoted to a very important technique in logarithmic geometry, which was introduced by Olsson in \cite{Olsson-logarithmic} (with a strong influence of ideas of Luc Illusie). It turns out that $S$-logarithmic structures on schemes $T$ over (the underlying scheme of) a base log scheme $S$ are classified by a stack $\Log_S$. Working with such stacks allows to interpret various logarithmic constructions and notions in terms of usual algebraic geometry of schemes and stacks. In particular, some results can be deduced from the non-logarithmic analogs on the nose.

\subsection{Constructions of $\Log_S$}

\subsubsection{The moduli definition of $\Log_S$} To any log scheme $S=(\uS,\calM_S)$ Olsson assigns the category $\Log_S$ fibered in groupoids over the category of $\uS$-schemes as follows: an object of $\Log_S$ is a logarithmic $S$-scheme $X$ and a morphism is a {\em strict} morphism $Y\to X$ of logarithmic $S$-schemes. The fiber functor just forgets the log structures.

\begin{rem}\label{Logrem}
(i) If $T$ is an $\uS$-scheme, then an object of $\Log_S(T)$ is just a logarithmic $S$-scheme whose underlying scheme is $T$. Thus, the stack $\Log_S$ parameterizes the ways in which one can enhance $\uS$-schemes with the structure of logarithmic $S$-schemes. So, informally speaking, it parameterizes log structures over $S$.

(ii) The association $S\mapsto\Log_S$ is naturally a functor from the category of schemes to the category of stacks.
\end{rem}

\subsubsection{Algebraicity}
Olsson proved that the stack $\Log_S$ is in fact an Artin stack of locally finite type over $S$. The proof goes by checking the usual properties -- representability of the diagonal and existence of a smooth presentation. The first property reduces to a simple study of the group of $\calM_S$-automorphisms of the log-structures log $S$-schemes -- naturally, they are extensions of diagonalizable groups by finite groups. The second property holds because $\coprod_{P\to Q}S_P[Q]\to\Log_S$ is smooth and surjective. The latter will be discussed in \S\ref{pressec} and then we will use it to construct $\Log_S$ very explicitly (in particular, the map from each $S_P[Q]$ factors through the quotient by the group of $S$-automorphisms, which is the extension of the diagonalizable group $\bfD_{Q^\gp/P^\gp}$ by the finite group $\Aut_P(Q)$).

\subsubsection{The tautological log structure}
By the definition of the stack $\Log_S$ any scheme over it is provided with a canonical log structure, and by descent one immediately obtains that the same is true for stacks over $\Log_S$. In particular, $\Log_S$ itself is provided with a tautological log structure $\calM$ and for any $x\in\Log_S(T)$ the induced homomorphism $(T,\calM_x)\to(\Log_S,\calM)$ is strict, that is, the log structure $\calM_x$ is induced from the tautological structure $\calM$ via the structure morphism $T\to\uS$. Similarly to Remark~\ref{Logrem} this provides a formalization of the claim that $\calM$ is the universal log structure over $S$. The following exercise essentially reduces to unravelling the definitions.

\begin{exer}
(i) Show that the log structure on $S$ induces a section $S\into\Log_S$ of the structure morphism $\Log_S\to S$.

(ii) Let $f\:Y\to X$ be an $S$-morphisms of logarithmic $S$-schemes and let $x\:X\to\Log_S$ and $y\:Y\to\Log_S$ be the corresponding morphisms. Show that $y=x\circ f$ if and only if $f$ is strict.

(iii) For a morphism $f\:Y\to X$ of log schemes the square
$$
\xymatrix{
\Log_Y\ar[d]\ar[r]^{\Log_f} & \Log_X\ar[d]\\
Y\ar[r]^f & X}
$$
is Cartesian if and only if $f$ is strict.
\end{exer}

\subsubsection{The stacks $\calX_P[Q]$}
The stack $\Log_X$ is huge, but it is a rather simple object that can be described by charts very explicitly. A first approximation for this is the following construction due to Olsson. Assume that $X$ is a log scheme with a global chart $X\to\bfA_P$ and $P\to Q$ is a homomorphism of monoids. Note that the diagonalizable group $\bfD_{Q^\gp/P^\gp}=\Spec(\bbZ[Q^\gp/P^\gp])$ acts on the $\bfA_P$-scheme $\bfA_Q$, and hence also acts on the $X$-scheme $X_P[Q]$ obtained by the base change. Let $\calX_P[Q]$ denote the quotient stack $[X_P[Q]/\bfD_{Q^\gp/P^\gp}]$.
The importance of these stacks introduced by Olsson becomes clear from the universal property they satisfy, which we are going to establish now. A short proof can be found in \cite[Lemma~2.2.4]{Molcho-Temkin}.

\begin{exer}
Let $X\to\bfA_P$ be a chart as above and let $Y$ be a logarithmic $X$-scheme.

(i) Show that flat locally on $Y$ the $P$-homomorphism $Q\to\ocalM_Y$ lifts to a homomorphism $Q\to\calM_Y$ and deduce that the functor of such liftings is a $D_{Q^\gp/P^\gp}$-torsor in the flat topology. Moreover, if $Y$ is fs, then this functor is already an \'etale torsor. (Hint: use that the homomorphisms $\calM^\gp_{Y}\onto\ocalM^\gp_{Y}$ has a section flat locally, and this is even true \'etale locally if $Y$ is fs and hence the groups $\ocalM^\gp_{Y,\oy}$ are torsion free.)

(ii) Deduce that $\calX_P[Q]$ represents the functor $\Hom_P(Q,\ocalM_Y)$ on the category of log schemes over $X$, while $X_P[Q]$ represents the functor $\Hom_P(Q,\calM_Y)$. In particular, $X$-homomorphisms $Y\to\calX_P[Q]$ are in a natural one-to-one correspondence with $P$-homomorphisms $Q\to\ocalM_Y$. (Hint: the second claim is clear, to deduce the first one divide by the action of $D_{Q^\gp/P^\gp}=\calX_P[Q^\gp]$ and use (i).)
\end{exer}

\begin{rem}
Keep the above notation and assume that $Q\to\calM_Y$ is a chart that induces an isomorphism $\phi\:Q=\Gamma(Y,\ocalM_Y)$. There are many liftings of $\phi$ to a chart of $Y$, obtained by multiplying monomials by units, but this is precisely the ambiguity which is killed by dividing by $\bfD_{Q^\gp/P^\gp}$. Thus $Y\to \calX_P[Q]$ can be viewed as the canonical $X$-chart of $Y$ determined only by the isomorphism $\phi$. Its only ambiguity is the group $\Aut_P(Q)$ of $P$-automorphisms of $Q$.
\end{rem}

\subsubsection{A smooth presentation}\label{pressec}
Assume that $X$ possesses a global chart $P\to\calM_X$. Olsson proves that the natural morphism $\coprod_{P\to Q}\calX_P[Q]\to\Log_X$, where the union is over all homomorphisms from $P$ to a fine monoid, is strict, surjective and \'etale. In particular, $\coprod_{P\to Q}X_P[Q]\to\Log_X$ is a smooth presentation of $\Log_X$. In fact, this easily reduces to the fact that \'etale locally any logarithmic $X$-scheme $Y$ possesses an \'etale cover $Y'\to Y$ whose source possesses a chart $Y'\to X_P[Q]$ and hence a strict morphism $Y'\to\calX_P[Q]$.

In general, $X$ possesses a global chart \'etale-locally, and the above construction is compatible with strict \'etale morphisms $X'\to X$. So, a presentation of $\Log_X$ can be obtained from a presentation of its fine enough strict \'etale cover.

\subsubsection{A groupoid presentation}\label{logpressec}
The presentation $f\:\coprod_{P\to Q}\calX_P[Q]\to\Log_X$ already gives a non-bad approximation of the source, but clearly it factors through $f'\:\coprod_{P\to Q}\calX_P[Q]/\Aut_P(Q)$. Even the \'etale morphism $f'$ is still not a monomorphism because $P$-automorphisms of localizations of $Q$ not necessarily come from $\Aut_P(Q)$, but it is easy to pin down the ambiguity -- one needs to identify all localizations of $Q$ in all possible ways, in particular, dividing $\calX_P[Q]$ by $\Aut_P(Q)$. Informally speaking, $\Log_S$ is obtained from the union of all charts $\coprod_{P\to Q}X_P[Q]$ by identifying all isomorphic open subcharts, in particular, dividing by automorphism: at first step this involves dividing by the groups $\bfD_{Q^\gp/P^\gp}$, and then by identifying all localizations of $Q_i$ and $Q_j$. In a sense, $\Log_X$ is nothing else but the universal $X$-chart constructed purely combinatorially.

Now let us outline the construction. It is more convenient to work geometrically, when the contravariant functor $Q\mapsto\calX_P[Q]$ is replaced by the functor $\Spec(Q)\to\calX_P[Q]$ from the category of affine Kato fans over $\Spec(P)$ because the latter globalizes in the obvious way. Moreover, one can naturally define a wider category of Kato stacks, and this functor extends to Kato stacks by (an appropriate) descent. Consider the diagram of all affine $P$-fans $\Spec(Q)$ with the morphisms being face embeddings, then the colimit $\calL_P$ exists as a Kato fan. Intuitively, it is a $P$-``fan'' which contains each $\Spec(Q)$ as a face in a unique way. It is not so difficult to show that $\Log_X=\calX_P[\calL_P]$, in particular, the morphism $\Log_X\to X$ is monoidal (in the stacky sense). Moreover, one can now give an explicit stacky presentation of $\Log_X$. We outline the main results in a (difficult) exercise below and refer to \cite[Sections 2,3]{Molcho-Temkin} for detailed arguments.

\begin{exer}
(i) Given $P$-monoids $Q_0\..Q_n$ by a join face we mean a $P$-monoid $R$ and face embeddings $\Spec(R)\into\Spec(Q_i)$ with $0\le i\le n$ (in other words, we fix isomorphisms of $R$ and localizations of $Q_i$). Show that the colimit $J(Q_0\..Q_n)$ of the diagram of all face embeddings of $Q_0\.. Q_n$ is a Kato fan, which we call the join of $Q_0\..Q_n$.

(ii) Construct a natural simplicial Kato fan $\calL_n=\coprod_{Q_0\.. Q_n}J(Q_0\..Q_n)$ and show that it is in fact a groupoid equivalent to a Kato fan which we denote $\calL_P$. In other words, $\calL_1=\coprod_{P\to Q}\Spec(Q)\to\calL_P$ is a cover and its fiber powers are $\calL_n$.

(iii) Show that $\calL_P$ is characterized by the following universal property: any $P$-fan $\Spec(Q)$ possesses a unique face embedding into $\calL$.

(iv) Show that $\calX_P[\calL_P]=\Log_X$ and deduce that $\Log_X$ is equivalent to the simplicial stack $\calX_n=\coprod_{Q_0\..Q_n}\calX_P[J(Q_0\..Q_n)]$, which is, in fact, a groupoid.
\end{exer}

\subsection{Stacks $\Log$ and logarithmic properties}
Now, we will show how one can systematically interpret various logarithmic properties of morphisms of log schemes. Initially such notions as log smoothness, log flatness, log \'etaleness, etc. were defined in a rather ad hoc manner. Then in \cite{Olsson-logarithmic} Olsson found a very general way to unify these definitions.

\begin{defin}
Let $\calP$ be a property of morphisms of schemes, for example, smooth, \'etale, flat, etc. A morphism of log schemes $f\:Y\to X$ is said to be {\em log $\calP$} (resp. {\em weakly log $\calP$}) if the associated morphism of stacks $\Log_f\:\Log_Y\to\Log_X$ (resp. $Y\to\Log_X$) is $\calP$.
\end{defin}

\begin{rem}
(i) Both definitions have some advantages. The morphism $Y\to\Log_X$ is ``smaller'' and easier to analyze; if $Y$ is quasi-compact, then a morphism $Y\to X$ factors through an open substack of $\Log_X$ finitely presented over $X$. On the other hand, when studying compositions it is certainly easier to work with the morphisms $\Log_Y\to\Log_X$.

(ii) Despite the terminology, neither condition implies the other one. Olsson showed in \cite[Example~4.3]{Olsson-logarithmic} that if $\calP$ is ``having geometrically connected fibers'', then log $\calP$ does not imply weakly log $\calP$, but there is even a much more basic example: any morphism $f\:Y\to X$ with a quasi-compact source is not weakly log surjective because $\Log_X$ is never quasi-compact, while $f$ is often log surjective, for example, when it is an isomorphism.

However, using that $Y\to\Log_X$ factors into the composition of an open immersion $Y\into\Log_Y$ and $\Log_f$ we obtain that if $\calP$ is local on the source, then log $\calP$ implies weakly log $\calP$.

(iii) For the properties of being log smooth, log \'etale and log flat Olsson showed that the original definitions given by Kato are equivalent to the new ones and also equivalent to the corresponding weak logarithmic properties.
\end{rem}

\begin{exer}
Use Olsson's definition to reprove Kato's criterion of log smoothness from Theorem~\ref{smoothchart}. Most probably, this will lead you to a conceptual explanation of the fact that condition \ref{smoothchart}(ii)(b) on $\phi\:P\to Q$ is equivalent to smoothness of the morphism $\bfD_{k(\oy),Q}\to\bfD_{k(\oy),P}$.
\end{exer}

Now we can also naturally interpret the notion of a log \'etale cover:

\begin{exer}
Prove that a morphism $f\:Y\to X$ is a log \'etale cover if and only if $\Log_f\:\Log_Y\to\Log_X$ is an \'etale cover. In other words, a log \'etale morphism is a cover if and only if it is log surjective.
\end{exer}

\subsubsection{Equivalence of the conditions}
A general result about equivalence of the two definitions was obtained in \cite[Theorem~4.3.1]{Molcho-Temkin}: if $\calP$ is stable under pullbacks, \'etale local on the source, and flat local on the base, then a morphism is log $\calP$ if and only if it is weakly log $\calP$. This follows from a slightly surprising fact that $\Log_Y\to\Log_X$ can be obtained from its small piece $Y\to\Log_X$ by base change and flat descent.

\begin{exer}
(i) Assume that $Y\to X$ has a global chart $\bfA_Q\to\bfA_P$. Show that for any homomorphism $Q\to R$ both squares in the following diagram are Cartesian:
\begin{align*}
\xymatrix{\calY_Q[R] \ar[d] & Y_Q[R] \ar[r] \ar[l] \ar[d] & Y \ar[d] \\ \calX_P[R] & [X_P[R]/\bfD_{Q^\gp/P^\gp}] \ar[r] \ar[l] & \calX_P[Q] }
\end{align*}

(ii) Using the presentation of the stacks $\Log$ constructed in \S\ref{logpressec}, \'etale descent and claim (i) deduce the assertion of \cite[Theorem~4.3.1]{Molcho-Temkin}.
\end{exer}

\subsubsection{Log regularity}
The above equivalence result applies, in particular, to the following properties: smoothness, \'etaleness, flatness and regularity. In fact, log regularity was introduced in \cite{Molcho-Temkin} and the equivalence with weak log regularity was used to establish its basic properties, e.g. a chart criterion analogous to Kato's chart criterion of log smoothness. In fact, this was the original motivation of the research of \cite{Molcho-Temkin}, which, in its turn, was motivated by relative resolution of singularities.

\subsubsection{Logarithmic differentials}
The stacks $\Log$ also allow to interpret logarithmic derivations and differentials. In view of the fact that a morphism $f\:Y\to X$ is log smooth if and only if the associated morphism $h\:Y\to\Log_X$ is smooth, the following fact is very natural: $\Omega^1_{Y/X}=\Omega^1_{\uY/\Log_X}$ and $\Der_{Y/X}=\Der_{\uY/\Log_X}$ (see \cite[Lemma~2.4.4]{ATW-relative}).

\begin{exer}
(i) Assume that $Y\to X$ possesses a chart $\bfA_Q\to\bfA_P$. Prove that $\Der_{Y/X}=\Der_{\uY/\calX_P[Q]}$. (Hint: compare the first fundamental sequences of log derivations associated with $Y\to X_P[Q]\to X$ and of derivations associated with $Y\to X_P[Q]\to\calX_P[Q]$.)

(ii) Use the \'etale cover $\coprod_Q\calX_P[Q]\to\Log_X$ to deduce that, indeed, for any morphism of log schemes $Y\to X$ one has that $\Der_{Y/X}=\Der_{\uY/\Log_X}$ and hence also $\Omega^1_{Y/X}=\Omega^1_{\uY/\Log_X}$.
\end{exer}

\subsubsection{Logarithmic fibers}
Let $f\:Y\to X$ be a morphism of log schemes. By the {\em log fibers} of $f$ we mean the connected components of the fibers of the induced morphism $h\:Y\to\Log_X$. For example, if $f$ is log smooth or log regular, then the log fibers are smooth or regular, respectively. It is easy to compute the log fibers \'etale-locally: if $f$ has a chart modeled on $\bfA_Q\to\bfA_P$, then $h$ factors through the \'etale morphism $\calX_P[Q]\to\Log_X$ hence the log fibers are nothing else but the fibers of the stacky chart $Y\to\calX_P[Q]$.

Let us consider two general types of examples of opposite kind. If the homomorphisms $\ocalM_x\to\ocalM_y$ are injective, then the log fibers have the most natural geometric interpretation:

\begin{exer}
Assume that $f\:Y\to X$ has a chart modeled on $\bfA_Q\to\bfA_P$ with sharp $P$ and $Q$ (such a chart exists \'etale-locally at $y\in Y$ with $x=f(y)$ if $\ocalM_x\to\ocalM_y$ is injective). Show that the log fibers of $f$ are the connected components of the log strata of the fibers of $f$. In particular, log fibers of a toroidal variety over a field are just the connected components of its log stratification.
\end{exer}

Such a description certainly cannot work for log blowings up, which might have non-discrete fibers but are log \'etale.

\begin{exer}
Show that the log fibers of any log blowing up $Y\to X$ are nothing else but the points of $y$ (as one would expect in the case of a monomorphism).
\end{exer}

In general, one can locally factor a morphism into a composition of a sharp morphism and a log blowing up, so log fibers admit a sort of a mixed description, but we do not discuss this here and refer the interested reader to \cite[\S2.2]{ATW-relative}.


\bibliographystyle{amsalpha}

\bibliography{lognotes}

\end{document}